\documentclass[12pt]{article}
\usepackage{amsmath,amsfonts,amssymb,amsthm}
\theoremstyle{plain}
\newtheorem{theorem}{Theorem}[section]

\theoremstyle{definition}
\newtheorem{definition}[theorem]{Definition}

\newcommand{\R}{{\mathbb R}}

\DeclareMathOperator{\card}{card}

\newcommand{\PP}{\mathbb{P}} % for probabilities
\newcommand{\EE}{\mathbb{E}} % for expectation

\begin{document}  % ========== BEGIN DOCUMENT

\title{On a $\psi_1$ - norm estimate of sum of dependent random variables using simple random walk on graph}

\author{Susanna Spektor}

\date{}

\maketitle
\begin{abstract}
We obtained a $\psi_1$ estimate for the sum of Rademacher random variables under condition that they are  dependent.
%%depends on the dimension of a vector.
\end{abstract}

%\thispagestyle{empty}

%\begin{abstract}

%\end{abstract}

\setcounter{page}{1}

\section{Introduction}

Let $X_1, \ldots, X_N$ be a sequence of independent real valued random variables and let $\Sigma=\sum_{i=1}^NX_i$.
The estimate of moments of $\Sigma$, that is of the quantities $\|\Sigma\|_p=\EE\left(\Sigma^p\right)^{1/p}$, are often appear in many areas of mathematics. The growth of moments is closely related to the behavior of the tails of $\Sigma$.

Probabilists have been interested in the moments of sums of random variables since the early part of last century. Khinchine’s 1923 paper appears to make the ﬁrst signiﬁcant contribution to this problem \cite{Khintchine}. It provides inequalities for the moments of a sum of Rademacher random variables. In 1970, Rosenthal generalised Khinchine’s result to the case of positive or mean-zero random variables \cite{Rosenthal}. Further reﬁnements to these bounds have been made by Latala and Hitczenko, Montogomery-Smith and Oleszkiewiez in more recent times \cite{Latala, Hitczenko}. Nowadays, it appears that in the different applications of  mathematics, statistics,  computer science and engineering similar estimates for the case when random variables are not independent are important (see for example \cite{Doukhan, PS}).

Our aim in the present work is to find bound on the sum of random variables, $\Sigma=\sum_{i=1}^{2n}X_i$, in the case when $X_i=a_i \varepsilon_i$, where $a \in \R^{2n}$ and $\varepsilon_i, i=1, \ldots, 2n$,  under an additional assumption on the Rademacher random variables, namely
 \begin{align}\label{2}
S=\sum_{i=1}^{2n}\varepsilon_i=0.
\end{align}
To shorter notation, by $\EE_S$ we  denote an expectation with assumption (\ref{2}).

Recall, the Rademacher random variables satisfying the following condition: $\PP(\varepsilon_i=1)=\PP(\varepsilon_i=-1)=\displaystyle{\frac 12}$, for $i=1,\ldots, 2n$. As usual for $\varepsilon \in \{\pm 1\}^{2n}$ by $\varepsilon_1, \ldots, \varepsilon_{2n}$ we denote coordinates of $\varepsilon$.

\bigskip

 Consider the following set

\begin{align}\label{set}
\Omega=\left\{\varepsilon \in\{-1, 1\}^{2n}\, |\,  \sum_{i=1}^{2n}\varepsilon_i=0\right\}=\left\{\varepsilon\in\{-1, 1\}^{2n}\, |\,  \card\{i: \, \varepsilon_i=1\}=n\right\}.
\end{align}
 Thus, for $\varepsilon \in \Omega$ the sequence of its coordinates is a sequence of a weekly dependent Rademacher random variables.

 For set $\Omega$ we put into correspondence the group $\Pi_{2n}$ of all permutations of set $\{1,..., 2n\}$ as
\begin{align*}
\sigma \in \Pi_{2n}\longleftrightarrow A_{\sigma}=\left\{\varepsilon \in \Omega \,\, | \,\, \varepsilon_i=1 \,\, \mbox{if} \,\,\sigma(i)\leq n; \, \varepsilon_i=-1 \,\,\mbox{if}\,\,\sigma(i)>n\right\}.
\end{align*}
Define $f:\Pi_{2n}\longrightarrow\R$ by
\begin{align}\label{function}
f(\sigma):=\left|\sum_{i=1}^na_{\sigma(i)}-\sum_{i=n+1}^{2n}a_{\sigma(i)}\right|.
\end{align}
Note, that $\EE_S\left|\sum_{i=1}^na_i \varepsilon_i\right|^p=\EE|f|^p$. Thus, it is enough to estimate $p$-th moments of $f$.

In the present paper we obtained the following result.
\begin{theorem}\label{th1}
Let $f$ defined as above. Then, for $p\geq 2$,
$$
(\EE f^p)^{1/p}\leq \EE |f| + C p\|a\|_2.
$$
\end{theorem}
%$$
% \left(\EE_S\left|\sum_{i=1}^{2n}a_i\varepsilon_i\right|^p\right)^{1/p} \leq
% C \, \min\left\{ p, \sqrt{ p \, \ln (2n)} \right\} \, \|a\|_{2}  .
%$$
% bounds for $\psi_1$ and $\psi_2$ norms of $\left|\sum_{i=1}^{2n}a_i\varepsilon_i\right|^p$

\bigskip

The paper is organized as following. In the next section we provide the necessary known tools and definitions. In Section 3, we will establish bounds on $\psi_1$-norm.

%%L$\acute{\textmd{e}}$vy

\bigskip

\section{Preliminaries}
\label{prelim}

\subsection{Orlicz norms and $\psi_{\alpha}$-estimates.}
%In this subsection we give an equivalent definition of the moment-estimate.

\begin{definition}
An \textit{{Orlicz function}} is a convex, increasing function $\psi: [0, \infty)\longrightarrow [0, \infty]$, such that $\psi(0)=0$ and $\psi(x)\longrightarrow \infty$ as $x\longrightarrow \infty$.
\end{definition}
Classical examples of  Orlicz functions are
\begin{align}\label{f}
\varphi_p(x)=\frac{x^p}{p}, \quad p\geq 1, \forall x\geq 0
\end{align}
and
\begin{align}\label{p}
\psi_{\alpha}(x)=e^{x^{\alpha}}-1, \quad \alpha\geq 1, \forall x \geq 0.
\end{align}

\begin{definition}
Let $\psi$ be an Orlicz function. For any real random variable $X$ on a measurable space $(\Omega, \sigma, \mu)$, define its \textit{$L_{\psi}$-norm} by
\begin{align*}
\|X\|_{\psi}:=\inf\{c>0: \EE \, \psi\left(|X|/{c}\right)\leq 1\}.
\end{align*}
We say $X$ is $\psi$-variable if $\|X\|_{\psi}< \infty$.
\end{definition}

The following well-known theorem describes the behaviour of a random variable with
bounded $\psi_{\alpha}$-norm (see for example \cite{Chafai}).
\begin{theorem}\label{psiesti}
Let $X$ be real-valued random variable and $\alpha \geq 1$. The following assertions are equivalent:
\begin{enumerate}
\item There exists $K_1>0$, such that $\|X\|_{\psi_{\alpha}}\leq K_1$.
\item There exists $K_2>0$, such that for every $p\geq \alpha$,
$$
\left(\EE|X|^p\right)^{1/p}\leq K_2p^{1/{\alpha}}.
$$
\item There exists $K_3, K_3'>0$, such that for every $t>K_3'$,
$$
\PP(|X|\geq t)\leq \exp\left(-t^{\alpha}/ K_3^{\alpha}\right).
$$
Note, $K_2\leq 2eK_1, \, K_3\leq e K_2, \, K_3'\leq e^2K_2, \ K_1\leq 2 \max(K_2, K_3')$.

\item In the case $\alpha>1$, let $\beta$ be such that $\displaystyle{\frac{1}{\alpha}+\frac{1}{\beta}=1}$. There exist $K_4, K_4'>0$ such that for every $\lambda\geq 1/K_4'$,
    $$
    \EE \, \exp\,(\lambda|X|)\leq \exp\,(\lambda K_4)^{\beta}.
    $$
    Note, $K_4\leq K_1,\,  K_4'\leq K_1, \, K_3'\leq 2K_4^{\beta}/(K_4')^{\beta-1}$.
\end{enumerate}
\end{theorem}
The space $L_{\psi}(\Omega, \sigma, \mu)=\{X: \|X\|_{\psi}< \infty\}$ is the Orlicz space associated to $\psi$.
Note that the Orlicz space associated to function $\varphi_p$, defined by (\ref{f}), is the classical $L_p$-space.

\subsection{Simple random walk on graph.}

 Let $G(V, E)$ be a connected undirected graph, where $V$ stays for a set of vertices and $E$ is a set of edges. A \textit{simple random walk} is a sequence of vertices $v_0, v_1, \ldots, v_t$, where $v_i \sim v_{i+1}$ (that is $\{v_i, v_{i+1}\}\in E$) for $i=0,1, \ldots, t-1$. That is, given an initial vertex $v_0$,  select randomly an adjacent vertex $v_1$, and move to this neighbor. Then, select randomly a neighbor $v_2$ of $v_1$, and move to it, etc. The probability it moves from vertex $v_i$ to $v_{i+1}$ (assuming it sits at $v_i$) is given by
\begin{align}\label{12}
p(v_i, v_{i+1})=\left\{ \begin{array}{rcl}
\displaystyle{\frac{1}{deg(v_i)}},& \quad \mbox{if} \hspace{0.1cm} v_i \sim v_{i+1}\\
\\
0,& \quad \mbox{otherwise},
\end{array}\right.
\end{align}
where $deg(v_i)$ denotes the degree of vertex $v_i$. This is a walk using a transition probability matrix, $P=\left(p(v_i, v_{i+1})\right)_{v_i, v_{i+1}\in V}$. The transition probability (\ref{12}) has a reversible equilibrium probability distribution $\mu(v_{i})$. That is,
\[
\mu(v_i)p(v_i, v_{i+1})=\mu(v_{i+1})p(v_{i+1}, v_i)
\]
and $\mu(v_i)$ is proportional to $deg(v_i)$.

Let $I$ be the $V \times V$ identity matrix. The discrete Laplacian is the matrix $L=P-I$ with its eigenvalues $0<\lambda_1\leq \lambda_2\leq \ldots$, ordered in non-increasing order. The smallest eigenvalue, $\lambda_1>0$, is called the \textit{spectral gap} of the random walk.

For $f: V\longrightarrow \R$ define
\begin{align}\label{13}
|||f|||^2_{\infty}=\frac 12 \sup_{v_i \in V}\sum_{v_{i+1}\in V}|f(v_i)-f(v_{i+1})|^2p(v_i, v_{i+1}).
\end{align}
We will use the following concentration inequality (see \cite{AidaStrook} or \cite{Ledoux}):
\begin{theorem}\label{t.h.cons}
Assume that $(p, \mu)$ is reversible on the finite graph $G(V,E)$, and let $\lambda_1>0$ be the spectral gap. Then, if $|||f|||_{\infty}^2< \infty$, we have
\begin{align}\label{14}
\mu\left(f> \int f d\mu+t\right)\leq 3 \exp\left({\frac{-t\sqrt{\lambda_1}}{2|||f|||^2_{\infty}}}\right).
\end{align}
\end{theorem}

\bigskip

Let us now specialize to $V=\Pi_{2n}$, the group of all permutations $\sigma$ of the set $\{1, \ldots, 2n\}$, and to $E=\{(\sigma, \sigma \tau)  \, \, \mid \, \, \tau\, \mbox{ is a transposition on } \Pi_{2n}\}$.
The transition probability $p(\sigma, \sigma \tau)$ on $G=(\Pi_{2n}, E)$ is
\begin{align}\label{15}
p(\sigma, \sigma \tau)=\frac{2}{(2n)^2},
\end{align}
and reversible equilibrium distribution $\mu$ on $\Pi_{2n}$ is a unique invariant measure for $p$ (see for example \cite{Chatterjee} for these facts). Also, as proved in \cite{DiaconisShahshahani}, the spectral gap of the random transposition walk on $\Pi_{2n}$ is $\displaystyle{\lambda_1=\frac{2}{2n}= \frac 1n}$. Thus, the concentration inequality (\ref{14}) for simple random walk on $G(\Pi_{2n}, E)$ can be rewritten as
\begin{align}\label{16}
\mu(\{\sigma: f(\sigma)-\EE f\geq t\})\leq \exp{\left(\frac{-t}{2|||f|||^2_{\infty} \sqrt n}\right)}.
\end{align}

\section{Proof of Theorem \ref{th1}}
\label{old}

%In this section we obtain a $\psi_1$-estimate for $\left|\sum_{i=1}^{2n}a_i\varepsilon_i\right|$
%under assumption (\ref{2}).

We are going to use inequality (\ref{16}). We calculate first
\[
|||f|||^2_{\infty}=\frac 12 \sup_{\sigma \in \Pi_{2n}} \, \sum_{\tau:\sigma \tau \in \Pi_{2n}}|f(\sigma)-f(\sigma \tau)|^2p(\sigma, \sigma \tau),
\]
where $p(\sigma, \sigma \tau)$ is defined in (\ref{15}).

Consider $g(\sigma)=\sum_{i=1}^na_{\sigma(i)}-\sum_{i=n+1}^{2n}a_{\sigma(i)}$. Since $\tau(i,j)$ is a random transposition with $i, j$ chosen uniformly from the set $\{1, \ldots, 2n\}$, we obtain
\[
g(\sigma)-g(\sigma \tau)=2(a_i-a_j)h(i,j),
\]
where
\[
\displaystyle{h(i,j)=\left\{ \begin{array}{rcl}
1&,  \quad \mbox{if} \hspace{0.1cm} j\leq n<i\leq 2n\\
-1&,  \quad \mbox{if} \hspace{0.1cm} i\leq n<j\leq 2n\\
0&,  \hspace{1.7cm} \mbox{otherwise}.
\end{array}\right.}
\]

Thus, $|f(\sigma)-f(\sigma \tau)|^2=4(a_i-a_j)^2h^2(i,j)$. And we can calculate
\begin{align*}
|||f|||^2_{\infty}&=\frac{1}{n^2}\sum_{\tau(i,j)}(a_i-a_j)^2h^2(i,j)\notag\\
&=\frac{2}{n^2}\sum_{i=1}^n\sum_{j=n+1}^{2n}(a_i-a_j)^2h^2(i,j)\notag\\
&=\frac{2}{n^2}\left(n\|a\|_2^2-2\sum_{i=1}^n\sum_{j={n+1}}^{2n}a_ia_j\right)
\end{align*}
Since
\[
\displaystyle{-\sum_{i=1}^n\sum_{j={n+1}}^{2n}a_ia_j\leq \sum_{i=1}^n\sum_{j={n+1}}^{2n}\frac{a_i^2+a_j^2}{2}=\frac n2 \|a\|^2_2},
\]
the last equation can be bounded by
\begin{align}\label{18}
|||f|||^2_{\infty}\leq \frac{4}{n}\|a\|^2_2.
\end{align}

Now, using (\ref{16}), (\ref{18}) and an upper bound $\Gamma(x)=x^{x-1}$, for all $x\geq 1$ (see for example \cite{Anderson}), we obtain
\begin{align*}
\EE (f-\EE f)^p=\int_0^{\infty}\mu((f(\sigma)-\EE f)^p\geq t^p) dt^p&\leq 6p \int_0^{\infty} e^{-t/(4\|a\|_2)}t^{p-1} dt\\
&=6 \, p\,  4^p \Gamma(p)\|a\|^p_2 \leq 4^p \, 6 p^p\|a\|^p_2.
\end{align*}
Hence
\[
(\EE f^p)^{1/p}\leq \EE |f| + 24 p\|a\|_2.
%\sqrt{\EE f^2}+Cp\|a\|_2.
\]

\noindent\textbf{Remark:} Note that $\displaystyle{\EE|f|\leq \left(\EE |f|^2\right)^{1/2}}$, where $\EE |f|^2$ can be directly calculated (see \cite{SS}).

\end{document}